\documentclass[12pt]{article}
\usepackage{amsmath,amsfonts}
\def\be{\begin{equation}} 
\def\ee{\end{equation}}
\newcommand{\beq}{\begin{eqnarray}}
\newcommand{\eeq}{\end{eqnarray}}
\newcommand{\nbeq}{\begin{eqnarray*}}
\newcommand{\neeq}{\end{eqnarray*}}
\def\D{\displaystyle}
\begin{document}

\title{On characterizations based on regression of linear
combinations of record values}
\author{George P. Yanev \\ {\it The University of Texas - Pan American, Edinburg,
USA}\\ \\ M. Ahsanullah  \\
{\it Rider University, Lawrenceville, USA} } \date{} \maketitle

\begin{abstract}
We characterize the exponential distribution in terms of the
regression of a record value with two non-adjacent record values as
covariates. We also study characterizations based on the
regression of linear combinations of record values. \\

\noindent {\it AMS (2000) Subject Classification}. 62G30, 62E10.

\end{abstract}

\begin{center}
\section{Introduction and results}
\end{center}

There is a number of studies on
 characterizations of probability distributions by means of regression relations of
 one record value with one
 or two other record values as covariates. For a recent review paper on the subject we refer to
 Pakes (2004), see also Ahsanullah and Raqab (2006), Chapter 6. To formulate and discuss
 our results we need to introduce some notations. Let $X_1, X_2,
\ldots$ be independent copies of a random variable $X$ with
absolutely continuous (with respect to the Lebesque measure) distribution function $F(x)$. An observation
$X_j$ is called a (upper) record value if it exceeds all previous
observations, i.e., $X_j$ is a (upper) record if $X_j>X_i$ for all
$i<j$. If we define the record times sequence by $T_1=1$ and $
T_n=\min \{j:X_j>X_{T_{n-1}}, j>T_{n-1}\}$, for $n>1$, then the
corresponding record values are $R_n=X_{T_n}$, $n=1,2,\ldots$
Let $F(x)$ be the cumulative distribution function of an exponential distribution given by
\begin{equation}
\label{exp_type} F(x)=1-e^{\D -c(x-l_F)}, \qquad -\infty<l_F\le x,
\end{equation} where $c>0$ is an arbitrary constant. Let us mention that (\ref{exp_type}) with $l_F>0$
 appears, for example, in reliability studies where $l_F$ represents the guarantee time; that is,
 failure cannot occur before $l_F$ units of time have elapsed (see Barlow and Proschan (1996), p.13).

Bairamov et al. (2005) study characterizations of exponential and
related distributions in terms of the regression of $R_n$ with two adjacent
record values as covariates.
They prove that
$F(x)$ is exponential if and only if \be \label{BAP05}
E[h'(R_n)|R_{n-1}=u, R_{n+1}=v ]=\frac{h(v)-h(u)}{v-u}, \qquad l_F<u<v, \ee
where the function $h$ satisfies some regularity conditions. Let us
note that by the mean-value theorem, there exists at least one
number $\xi$ inside the interval $(u,v)$ such that for
$u>l_F$
\[ E[h'(R_n)|R_{n-1}=u, R_{n+1}=v ]=h'(\xi).\]
In particular, if $h'(x)=x$ then (\ref{BAP05}) becomes \[
E[R_n|R_{n-1}=u, R_{n+1}=v ]=\frac{u+v}{2}, \qquad l_F<u<v.\]

\noindent Yanev et al. (2008) extend (\ref{BAP05}) to the case when at least one of the
two covariates is adjacent to $R_n$.
To formulate their result, we need to introduce some notations. Further on, for a given $h(x)$, we
denote
\[ \label{M_not}
M(u,v)=\frac{h(v)-h(u)}{v-u}, \
_iM_j(u,v)=\frac{\partial^{i+j}}{\partial u^i\partial
v^j}\left(\frac{h(v)-h(u)}{v-u}\right),  u\ne v,
\] as well as $_iM(u,v)$ and $M_j(u,v)$ for the $i$th
and $j$th partial derivative of $M(u,v)$ with respect to $u$ and
$v$, respectively. Let $k$ and $n$ be integers, such that $1\le k\le n-1$. When at least one of the covariate record values is
adjacent to $R_n$, it is shown in \cite{YAB}, under some
regularity assumptions, that $F(x)$ is exponential if and only if  \be
\label{YAB08}E[h^{(k)}(R_n)|R_{n-k}=u, R_{n+1}=v]= kM_{k-1}(u, v), \qquad l_F<u<v.
\ee In particular, if $h(x)=x^{k+1}/(k+1)!$, then $h^{(k)}(x)=x$ and
(\ref{YAB08}) becomes  \be \label{weights}
E[R_n|R_{n-k}=u, R_{n+1}=v]= \frac{u+kv}{k+1}, \qquad l_F<u<v. \ee Observe that the
right-hand side of (\ref{weights}) is the weighted mean of the covariates with
weights equal to the number of spacings they are away from $R_n$. One can also
see (using the arguments in \cite{YAB}) that for $r\ge 1$ \be
\label{weights2} E[R_n|R_{n-1}=u, R_{n+r}=v]= \frac{ru+v}{r+1}, \qquad l_F<u<v \ee
characterizes the exponential distribution too.

 If both
covariates are non-adjacent to $R_n$ the situation is more
complex.  Let $k$, $r$, and $n$ be integers, such that $1\le k\le n-1$ and $r\ge 1$.   Yanev et al (2008) obtain a necessary condition for exponentiality of $F(x)$. Namely, they prove, under some regularity
assumptions,  that if $F(x)$ is exponential, then for $l_F<u<v$
\[
E[h^{(k+r+1)}(R_n)|R_{n-k}=u,
R_{n+r}=v]=\frac{(k+r-1)!}{(k-1)!(r-1)!}\  _{r-1}M_{k-1}(u,v).
\]
However, no sufficient condition for $F(x)$ to be exponential that
involves only single regression of $R_n$ on two non-adjacent covariates is known yet. For example, in
\cite{YAB} the necessary and sufficient condition for $F(x)$ to be
exponential is that both \[
E[R_n|R_{n-k}=u, R_{n+r}=v]=\frac{ ru+kv}{k+r}, \qquad l_F<u<v\]
and \[
 E[R_n|R_{n-k+1}=s, R_{n+r}=v]=\frac{rs+(k-1)v}{k+r-1}, \qquad l_F<s<v\] hold. Our first result provides new sufficient
 and necessary conditions for
(\ref{exp_type}) when both covariates are non-adjacent to $R_n$.
The conditions given
below are written in a form which extends (\ref{YAB08}). They
are alternative to and more compact than the results
in Theorem 1B of \cite{YAB}. We have the following theorem.

{\bf Theorem 1.}\ {\it Let $k$, $r$, and $n$ be integers, such that $1\le k\le n-1$ and $r\ge 1$. Assume that $F(x)$ is absolutely continuous.
Suppose $h(x)$ satisfies

(i) $h(x)$ is continuous in
$[l_F, \infty)$ and $h^{(k+r-1)}(x)$ is continuous in $(l_F,\infty)$;

(ii)
$_{r-1}M_k(l_F,v)\ne 0$ for $ v> l_F$.

 \noindent Then $F(x)$ is the
exponential cdf (\ref{exp_type}) if and only if for $l_F<u<s<v$ \beq
\label{new_res}
\lefteqn{\hspace{-1.5cm}\frac{k-1}{_{r-1}M_{k-1}(u,v)}E\left[h^{(k+r-1)}(R_n)
 \ | \ R_{n-k}=u, R_{n+r}=v\right]} \\
& = & \frac{k+r-1}{_{r-1}M'_{k-2}(s,v)}E\left[h^{(k+r-1)}(R_n)
\ | \ R_{n-k+1}=s, R_{n+r}=v\right], \nonumber\eeq where
$M'(u,v)=[h'(v)-h'(u)]/(v-u)$.}

Notice that, setting  $r=1$ and $k=2$ and letting $s\to v^-$,
one can see that (\ref{new_res}) reduces to (\ref{YAB08}) for the case $k=2$. We
illustrate the applicability of Theorem~1 with one corollary below.
Let $h(x)=x^{k+r}/(k+r)!$ and thus $h^{(k+r-1)}(x)=x$.
It is not difficult
  to see that with this choice of $h(x)$
  \[
  M(u,v)=\frac{v^{k+r-1}+\ldots +v^{k}u^{r-1}+v^{k-1}u^r+\ldots
+u^{k+r-1}}{(k+r)!}
\]
and
\[
\frac{(k+r-1)!}{(k-1)!(r-1)!}\
_{r-1}M_{k-1}(u,v)=\frac{ru+kv}{k+r}.
\]
Now, Theorem~1 implies the following corollary.

{\bf Corollary 1}\ {\it Let $k$, $r$, and $n$ be integers, such that $1\le k\le n-1$ and $r\ge 1$. Suppose $F(x)$ is absolutely continuous.
$F(x)$ is the exponential cdf (\ref{exp_type}) if and only if for $l_F<u<s<v$
\[
 \frac{k+r}{ru+kv}E\left[R_n | R_{n-k}=u, R_{n+r}=v\right]\!
=\! \frac{k+r-1}{rs+(k-1)v}E\left[R_n
|R_{n-k+1}=s, R_{n+r}=v\right].\]  }

Corollary~1 gives a characterization of (\ref{exp_type}), which is
alternative to that in Theorem 2B of \cite{YAB}, mentioned before Theorem~1
above.

Next we turn our attention to characterizations based on
regressions of differences (spacings) of two record values. Consider the
Weibull distribution given for $\alpha>0$ by its cdf \be
\label{Weibull}F(x)=1-e^{\D -cx^\alpha},  \qquad x\ge 0,
\end{equation} where $c>0$ is an arbitrary constant.
 Akhundov and Nevzorov (2008)
 study the regression of spacings of
record values as follows \be \label{AN} E[R_3-R_2|R_{1}=u,
R_{4}=v]=\frac{v-u}{3}, \qquad u<v. \ee If $F(x)$ is the exponential
(\ref{exp_type}) then it is clear that (\ref{weights}) and
(\ref{weights2}) lead to (\ref{AN}). Since (\ref{AN}) is a weaker condition
than (\ref{weights}) (or (\ref{weights2})), it is not a sufficient
condition for $F(x)$ to be exponential. Akhundov and Nevzorov
(2008) prove the interesting fact that there is only one more
 family of distributions, other than the exponential, that satisfies
(\ref{AN}). It turns out that (\ref{AN}) holds if and only if
$F(x)$ satisfies (\ref{Weibull}) with either $\alpha=1$ or
$\alpha=1/2$. Making use of the findings in \cite{YAB} and
Theorem~1 above, we generalize this result in two directions: (i)
considering $R_m-R_n$ for any $2\le m\le n-1$; and (ii) in the case of non-adjacent
covariates.
The following
characterization holds.

{\bf Theorem 2.}\ {\it Let $k$, $r$, $m$, and $n$ be integers, such that $1\le k\le m-1$, $r\ge 1$, and $2\le m\le n-1$. Suppose $F(x)$ is absolutely continuous in
$(0,\infty)$. Then $F(x)$ is given by (\ref{Weibull}) with
$\alpha=1$ or $\alpha=1/2$ if and only if for
$0<u<s<t<v<\infty$ \be \label{new_stuff1} \hspace{-3.5cm}
\frac{d+2}{v-u}\  E\left[
R_n-R_m |\ R_{m-k}=u, R_{n+r}=v\right]
\ee
\[
 =\frac{d}{t-s}\ E\left[ R_n-R_m| R_{m-k+1}=s,
R_{n+r-1}=t\right],
\]
where $d=n-m+k+r-2$.
}

Setting $n=m+1$ in (\ref{new_stuff1}) we obtain the following result.

{\bf Corollary 2.}\ {\it Let $k$, $r$, and  $m$ be integers, such that $1\le k\le m-1$, $r\ge 1$. Suppose $F(x)$ is absolutely continuous in
$(0,\infty)$. Then $F(x)$ is given by
(\ref{Weibull})
 with
$\alpha=1$ or $\alpha=1/2$ if and only if for $0<u<v<\infty$ \[
E\left[R_{m+1}-R_{m}\ |\
R_{m-k}=u, R_{m+r+1}=v\right]= \frac{k+r-1}{k+r+1}(v-u).\]}

Setting $k=r=1$ in (\ref{new_stuff1}) we obtain a corollary for
adjacent covariates.

{\bf Corollary 3.}\ {\it Let $m$ and $n$ be integers, such that $2\le m\le n-1$. Suppose $F(x)$ is absolutely continuous in
$(0,\infty)$. Then $F(x)$ is given by
(\ref{Weibull})
 with
$\alpha=1$ or $\alpha=1/2$ if and only if for $0<u<v<\infty$ \[
E\left[R_n-R_m\ |\
R_{m-1}=u, R_{n+1}=v\right]= \frac{n-m}{n-m+2}(v-u).\]}

Corollary 3 can be interpreted as follows. Let fix the integers $m$ and $n$, such that $2\le m<n$. According to
(\ref{weights}) and (\ref{weights2}), the conditional expectations of the spacings $R_m-R_{m-1}$ and $R_{n+1}-R_n$ given
$R_{m-1}=u$ and $R_{n+1}=v$  are equal and their sum is $2(v-u)/(n-m+2)$, that is
\nbeq
E[R_m-R_{m-1}|R_{m-1}=u,R_{n+1}=v] & = &
    E[R_{n+1}-R_n|R_{m-1}=u, R_{n+1}=v] \\
    & = & \frac{v-u}{n-m+2},
    \neeq
if and only if  $F(x)$ is exponential. Now, assume that it is only known that the conditional expectations above
have a sum $2(v-u)/(n-m+2)$, that is
\beq  \label{sum_spec}
\lefteqn{\hspace{-5cm}E[R_m-R_{m-1}|R_{m-1}=u,R_{n+1}=v] +
    E[R_{n+1}-R_n|R_{m-1}=u, R_{n+1}=v]}\nonumber\\
    & & =\frac{2(v-u)}{n-m+2}.
    \eeq
This is equivalent to
\nbeq
\lefteqn{E[R_n-R_m|R_{m-1}=u,R_{n+1}=v]}\\
& = &
        E[R_{n+1}-R_{m-1}|R_{m-1}=u,R_{n+1}=v]- E[R_m-R_{m-1}|R_{m-1}=u,R_{n+1}=v]\\
        & & - E[R_{n+1}-R_n|R_{m-1}=u,R_{n+1}=v]\\
        & = & v-u - \frac{2(v-u)}{n-m+2}\\
        & = & \frac{n-m}{n-m+2}(v-u).
        \neeq
        Therefore, according to Corollary~3, (\ref{sum_spec}) holds if and only if the underlying distribution is either exponential or Weibull  with $\alpha=1/2$.

Finally, we investigate the regression
\[
\frac{1}{cv-du}E\left[ aR_n-bR_m\ |\ R_{m-k}=u,
R_{n+r}=v\right], \qquad l_F<u<v,
\]
where  $a$, $b$, $c$, and $d$ with $a\neq b$ are some real
numbers. What choice of these numbers characterizes the
exponential distribution alone? The theorem below answers this
question.

 {\bf Theorem 3.} {\it Let $k$, $r$, $m$, and $n$ be integers, such that $1\le k\le m-1$, $r\ge 1$, and $2\le m\le n-1$. Suppose $F(x)$ is
absolutely continuous. Then $F(x)$ is exponential given by
(\ref{exp_type}) if and only if   for $l_F < s<u<v<t$ \beq
\label{newstuff2}
\lefteqn{\hspace{-1cm}E\left[ kR_n-(n-m+k)R_m|R_{m-k}=u, R_{n+r}=v\right]}\\
& & \hspace{-0.7cm} =\!\frac{du}{(d+1)s-t}E\left[
kR_n-(n-m+k)R_m |R_{m-k+1}=s,
R_{n+r-1}\!= \!t\right],  \nonumber \eeq where   $d=n-m+k+r-2$.  }

Setting $k=r=1$ in (\ref{newstuff2}) we obtain
\[
E\left[R_n-(n-m+1)R_m\ |\ R_{m-1}=u, R_{n+1}=v \right] =(m-n)u.
\]
Therefore, we have the following corollary.

 {\bf Corollary 4.} {\it Let $m$ and $n$ be integers, such that $2\le m\le n-1$. Suppose $F(x)$ is
absolutely continuous. Then $F(x)$ is exponential given by
(\ref{exp_type}) if and only if for $l_F<u<v$
\beq \label{cor4} \lefteqn{ E\left[R_n-(n-m)R_m\ |\ R_{m-1}=u, R_{n+1}=v\right]}\\
& & =E\left[R_m-(n-m)R_{m-1}|R_{m-1}=u,R_{n+1}=v\right]. \nonumber
\eeq }
It is interesting to note that setting $n=m+1$ the condition (\ref{cor4}) becomes
\[
E[R_{m+1}-R_m|R_{m-1}=u,R_{m+2}=v]=E[R_m-R_{m-1}|R_{m-1}=u,R_{m+2}=v].
\]

We shall prove the results presented here in the next three sections.

\setcounter{equation}{0}

\hspace{0.5cm} {\centering \section{Proof of Theorem 1}}

{\bf Sufficiency}.\ Denote the cumulative hazard function of the cdf $F(x)$ by $H_x=-\ln(1-F(x))$. Also, for simplicity, we will sometimes write
$W_{x,y}=H_y-H_x$. Referring to the Markov
dependence of the record values, one can show (e.g., Ahsanullah (2008)) that
the conditional density of $R_n$ given $R_{n-i}=u$
 and $R_{n+j}=v$  is  for $1\leq i\le n-1$ and $j\geq 1$
 \begin{equation}
\label{cond_den} \frac{\displaystyle (i+j-1)!}{\displaystyle
(i-1)!(j-1)!}\frac{\displaystyle W_{u,t}W_{t,v}}{\displaystyle
W_{u,v}^{i+j-1}}
H'_t.
\end{equation}
Assuming (\ref{new_res}), we will show that $F(x)$ satisfies
(\ref{exp_type}). Denote $d=k+r-2$. Referring to (\ref{cond_den}), it is not
difficult to obtain \nbeq \lefteqn{E\left[h^{(d+1)}(R_n)
 |R_{n-k}=u, R_{n+r}=v\right]}\\
 & = &
\frac{(d+1)!}{(k-1)!(r-1)!W_{u,v}^{d+1}} \int_u^v
h^{(d+1)}(x)W_{u,x}^{k-1}W_{x,v}^{r-1}dH_x  \\
& = &   \frac{(d+1)!}{(k-1)!(r-1)!W_{u,v}^{d+1}}I(u,v; k, r)
\quad \mbox{say},  \neeq and
\nbeq\lefteqn{E\left[h^{(d+1)}(R_n)
 |R_{n-k+1}=s, R_{n+r}=v\right]}\\
 & = &
\frac{d!}{(k-2)!(r-1)!W_{s,t}^d} \int_s^v
h^{(d+1)}(x)W_{s,x}^{k-2}W_{x,v}^{r-1}dH_x \\
& = &  \frac{d!}{(k-2)!(r-1)!W_{s,t}^d}I(s,v; k-1,r)
\quad \mbox{say}.\neeq
Now, one can see that (\ref{new_res}) is equivalent to \be
\label{new_mediate}\frac{W_{s,v}^d
\ _{r-1}M'_{k-2}(s,v)}{I(s,v;k-1,r)}I(u,v;k,r)
 =\ _{r-1}M_{k-1}(u,v)W_{u,v}^{d+1}. \ee
Differentiating (\ref{new_mediate}) with respect to
$u$ and letting $s\to u^+$, we have \nbeq \lefteqn{-(k-1)\
_{r-1}M'_{k-2}(u,v)W_{u,v}^dH'_u}\\
& & =\ _rM_{k-1}(u,v)W_{u,v}^{d+1}-(k+r-1)\
_{r-1}M_{k-1}(u,v)W_{u,v}^dH'_u. \neeq Dividing by
$W_{u,v}^{d+1}$ and grouping, we arrive at the equation \be
\label{eqn}\frac{H'_u}{H_u-H_v}=\frac{\ _rM_{k-1}(u,v)}{(k-1)\
_{r-1}M'_{k-2}(u,v)-(k+r-1)\ _{r-1}M_{k-1}(u,v)}, \ee provided that
the denominator in the right-hand side is not zero. (This is
equivalent to the assumption $\ _rM_{k-1}(u,v)\ne 0$, as we will see below.)
Since
 \begin{eqnarray*}
_{r-1}M'_{k-2}(u,v) & = & \ \frac{\partial^{k+r-3}}{\partial
u^{r-1}\partial v^{k-2}}\left[\frac{h'(v)-h'(u)}{v-u}\right]\\
    & = &
 \frac{\partial^{k+r-3}}{\partial
u^{r-1}\partial v^{k-2}} \left[M_1(u,v)+\ _1M(u,v)\right]\\
    & = &
\ _{r-1}M_{k-1}(u,v)+\ _{r}M_{k-2}(u,v), \end{eqnarray*} for the
denominator in the right-hand side of (\ref{eqn}) we have
\begin{eqnarray} \lefteqn{\hspace{-0.5cm}(k-1)\
_{r-1}M'_{k-2}(u,v)-(k+r-1)\
_{r-1}M_{k-1}(u,v)}\label{denominator}\\
 & = &
(k-1)[\ _{r-1}M_{k-1}(u,v)+\ _{r}M_{k-2}(u,v)]-(k+r-1)\
_{r-1}M_{k-1}(u,v)
 \nonumber \\
    & = &
    (k-1)\ _{r}M_{k-2}(u,v)-r\
    _{r-1}M_{k-1}(u,v) \nonumber \\
    & = &
    \ _{r}M_{k-1}(u,v)(u-v). \nonumber
\end{eqnarray} The last equality follows from Lemma~1 in \cite{YAB}. Now, (\ref{eqn}) and (\ref{denominator}) imply
\begin{eqnarray*} \frac{\displaystyle
H'_u}{\displaystyle H_u-H_v}
    & = &
    \frac{\displaystyle 1}{\displaystyle u-v}.
    \end{eqnarray*}
Integrating both sides with respect to $u$ from $l_F$ to $v$, we
obtain
\[
\ln(H_v-H_{l_F})=\ln (v-l_F) + \ln c,    \qquad c>0
\]
and thus $ H_v=c(v-l_F)$ which implies (\ref{exp_type}).

{\bf Necessity}. According to Theorem 1B in \cite{YAB}, if $F(x)$
satisfies (\ref{exp_type}), then 
\[ E[h^{(k+r-1)}(R_n)|R_{n-k}=u, R_{n+r}=v]=
\frac{(k+r-1)!}{(k-1)!(r-1)!}\ _{r-1}M_{k-1}(u, v)\]
and 
\[ E[h^{(k+r-1)}(R_n)|R_{n-k+1}=s, R_{n+r}=v]=
\frac{(k+r-2)!}{(k-2)!(r-1)!}\ _{r-1}M'_{k-2}(s, v). \]  These two
equalities imply (\ref{new_res}). This completes the proof of the
theorem.

\setcounter{equation}{0}

\hspace{0.5cm} {\centering \section{Proof of Theorem 2}}

To prove Theorem 2 we will need the following three lemmas.

{\bf Lemma 1} (Akhundov and Nevzorov (2008)) {\it Let $F(x)$ be
absolutely continuous. The equation
\[
H_v-H_u=\frac{2H'_uH'_v}{H'_u+H'_v}(v-u),\qquad 0<u<v<\infty
\]
has exactly two solutions given  by
$F(x)=1-\exp\{-cx^\alpha\}$ for $\alpha=1$ or $\alpha=1/2$, where $c>0$ is an
arbitrary constant. }

The following lemma is a straightforward corollary of Lemma~2 in
\cite{YAB}.

{\bf Lemma 2}\ {\it Let $k$, $r$, and $n$ be integers such that $1\le k\le n-1$ and $r\ge 1$. If $ F(x)=1-\exp\{-c(x-l_F)\}$,
$(l_F<x<\infty)$, where $c>0$ is an arbitrary constant, then
\[
E[R_n|R_{n-k}=u, R_{n+r}=v]=\frac{ru+kv}{k+r}, \qquad l_F<u<v<\infty.
\]
}

{\bf Lemma 3}\ {\it Let $a$ and $b>a$ be real numbers
and $i$ and $j$ be positive integers. Then
\[
I=\int_a^b [(y-a)^j(b-y)^i-(y-a)^i(b-y)^j]y^2\;
dy=\frac{i!j!(j-i)}{(i+j+2)!}(b-a)^{i+j+1}(b^2-a^2).
\]}

{\bf Proof.}\  We have
\beq \label{integral}
I & = & \int_a^b (y-a)^j(b-y)^iy^2\; dy - \int_a^b(y-a)^i(b-y)^jy^2\; dy \\
    & = &
    I_1-I_2 .\nonumber
    \eeq
Making in $I_1$ the change of variables $w=(y-a)/(b-a)$, we obtain
\be \label{first} I_1=(b-a)^{i+j+1}\int_0^1
w^j(1-w)^i[(b-a)w+a]^2dw. \ee Similarly, making in $I_2$ the change
of variables $w=(b-y)/(b-a)$, we have \be \label{second}
I_2=(b-a)^{i+j+1}\int_0^1 w^j(1-w)^i[b-(b-a)w]^2dw. \ee From
(\ref{first}) and (\ref{second}), we have \nbeq
I   & = &
    (b-a)^{i+j+1}\int_0^1w^j(1-w)^i\{ [(b-a)w+a]^2-[b-(b-a)w]^2\}dw\\
    & = &
    (b-a)^{i+j+2}(b+a)[2B(j+2, i+1)-B(j+1,i+1)]\\
    & = &
    (b-a)^{i+j+1}(b^2-a^2)B(j+1,i+1)\left(2\frac{j+1}{i+j+2}-1\right)\\
    & = &
(b-a)^{i+j+1}(b^2-a^2)\frac{i!j!(j-i)}{(i+j+2)!},
\neeq
which proves the lemma.

\vspace{0.5cm}{\bf Proof of Theorem 2.}\ {\bf Sufficiency.} We
shall  prove that (\ref{new_stuff1}) implies (\ref{Weibull}) with
either $\alpha=1$ or $\alpha=1/2$. First, assume that
$1\le k\le m-1$ and $r \ge 2$.
Referring to (\ref{cond_den}), one can obtain (recall that $d=n-m+k+r-2$) \beq \label{tope}
\lefteqn{E[R_n-R_m|R_{m-k}=u, R_{n+r}=v]}\\
& = & \frac{(d+1)!}{W_{u,v}^{d+1}}
\int_u^v \left[
\frac{W_{u,x}^{d-r+1}W_{x,v}^{r-1}}{(d-r+1)!(r-1)!}
-\frac{W_{u,x}^{k-1}W_{x,v}^{d-k+1}}{(k-1)!(d-k+1)!}\right]xdH_x
 \nonumber \\
& = & \frac{(d+1)!}{W_{u,v}^{d+1}}\ I(u,v; k,
r),\quad \mbox{say}, \nonumber \eeq and
\beq \label{bote} \lefteqn{E[R_n-R_m|R_{m-k+1}=u, R_{n+r-1}=v]}\\
& = &  \frac{(d-1)!}{W_{s,t}^{d-1}}
 \int_s^t \left[
\frac{W_{s,x}^{d-r}W_{x,t}^{r-2}}{(d-r)!(r-2)!}- \frac{W_{s,x}^{k-2}W_{x,t}^{d-k}}{(k-2)!(d-k)!}\right]xdH_x
\nonumber \\
& = &  \frac{(d-1)!}{W_{s,t}^{d-1}}\ I(s,t;
k-1,r-1),\quad \mbox{say}. \nonumber \eeq

Now, making use of (\ref{tope}) and (\ref{bote}), we can write
(\ref{new_stuff1}) as
 \be
\label{fraction2}
\frac{(d+1)!(t-s)W_{s,t}^{d-1}}{(d-1)!
I(s,t;k-1,r-1)}I(u,v;k,r)
= (v-u)W_{u,v}^{d+1}.\ee Let us
differentiate both sides of (\ref{fraction2}) with respect to $u$
and $v$. Then, after letting $s\to u^+$ and $t\to v^-$, (\ref{fraction2}) simplifies to
\[
d(v-u)H'_uH'_v=
(H'_u+H'_v)(H_v-H_u)
+d(v-u)H'_uH'_v. \]Therefore, \be \label{AN_eqn}
H_v-H_u=\frac{2H'_uH'_v}{H'_u+H'_v}(v-u). \ee According to Lemma~1,
equation (\ref{AN_eqn}) has the two solutions given by
(\ref{Weibull}) with $\alpha=1$ or $\alpha=1/2$. In the case $k=1$ and
$r\ge 2$ the proof is similar and is omitted here. If $k=r=1$, then
(\ref{new_stuff1}) simplifies to\[
\label{case_k=1} \frac{d+2}{d}E\left[
R_n-R_m|R_{m+1}=u,R_{n-1}=v\right]=v-u. \] Repeating the arguments for the case $k\ge 2$ above, it is not difficult to obtain equation (\ref{AN_eqn}).
The sufficiency is proved.

{\bf Necessity.} We need to show that both cdf's $F_1(x)=1-\exp\{-cx\}$
and $F_2(x)=1-\exp\{-cx^{1/2}\}$ satisfy (\ref{new_stuff1}). In
case of $F_1(x)$, it is not difficult to obtain the relation
(\ref{new_stuff1}) using Lemma~2 above. It remains to prove that
$F_2(x)=1-\exp\{-cx^{1/2}\}$ satisfies (\ref{new_stuff1}) as well. First assume that $2\le
k\le m-1$ and $r\ge 2$.
Since $x=H^2(x)/c^2$, for the left-hand side of (\ref{new_stuff1})
we have \nbeq
\lefteqn{\frac{d+2}{v-u} E\left[ R_n-R_m\ |\ R_{m-k}=u, R_{n+r}=v\right]}\\
& & \hspace{-0.8cm} = \frac{(d+2)!}{c^2W_{u,v}^{d+1}(v-u)}
\int_{H_u}^{H_v} \left[
\frac{W_{u,x}^{d-r+1}W_{x,v}^{r-1}}{(d-r+1)!(r-1)!}
-\frac{W_{u,x}^{k-1}W_{x,v}^{d-k+1}}{(d-k+1)!(k-1)!}\right]H_x^2dH_x.
\neeq \noindent Using Lemma~3 (twice) with $a=H_u$, $b=H_v$, after
some algebra, we obtain \be \label{LHS}
\frac{d+2}{v-u} E\left[ R_n-R_m\ |\ R_{m-k}=u, R_{n+r}=v\right]=2(n-m) . \ee

\noindent Similarly, using Lemma~3 with $a=H_s$ and  $b=H_t$ for
the right-hand side of (\ref{new_stuff1}) we have \be\label{RHS}
\frac{d}{t-s}E\left[ R_n-R_m\ |\ R_{m-k+1}=s, R_{n+r-1}=t\right]
= 2(n-m). \ee It follows from
(\ref{LHS}) and (\ref{RHS}) that $F_2(x)$ satisfies
(\ref{new_stuff1}). When $k=1$ and $r\ge 2$ or $k=r=1$ the proof is similar and
is omitted here.

\setcounter{equation}{0}

\hspace{0.5cm} {\centering \section{Proof of Theorem 3}}

 {\bf Sufficiency.} We shall prove that
(\ref{newstuff2}) implies (\ref{exp_type}).
First, assume that $2\le k\le
m-1$ and $r\ge 2$. Referring to (\ref{cond_den}) we obtain \nbeq  \label{thm31}
\lefteqn{E\left[ kR_n-(n-m+k)R_m\ | \ R_{m-k}=u, R_{n+r}=v\right]}\\
& & \hspace{-0.8cm}= \frac{(d+1)!}{W_{u,v}^{d+1}}
\int_u^v \left[
\frac{kW_{u,x}^{d-r+1}W_{x,v}^{r-1}}{(d-r+1)!(r-1)!}
-\frac{(d-r+2)W_{u,x}^{k-1}W_{x,v}^{d-k+1}}{(k-1)!(d-k+1)!}\right]xdH_x
 \nonumber \\
& &
\hspace{-0.8cm} = \frac{(d+1)!}{W_{u,v}^{d+1}}J(u,v;k,r),\quad
\mbox{say,  } \nonumber
 \neeq
 \noindent and
\nbeq \label{thm32} \lefteqn{E\left[
kR_n-(n-m+k)R_m \ |\ R_{m-k+1}=s,
R_{n+r-1}= t\right]}\\
& & \hspace{-0.8cm} = \frac{(d-1)!}{W_{s,t}^{d-1}}
\int_s^t
\left[
\frac{kW_{s,x}^{d-r}W_{x,t}^{r-2}}{(d-r)!(r-2)!}
-
\frac{(d-r+2)W_{s,x}^{k-2}W_{x,t}^{d-k}}{(k-2)!(d-k)!}\right]xdH_x
 \nonumber \\
& & \hspace{-0.8cm}=
\frac{(d-1)!}{W_{s,t}^{d-1}}J(s,t;k-1,r-1),\quad
\mbox{say}. \nonumber
 \neeq
Now, we can write (\ref{newstuff2}) as
 \be
\label{fraction3}
\frac{(d+1)!((d+1)s-t)W_{s,t}^{d-1}}{(d-1)!
J(s,t;k-1,r-1)}J(u,v;k,r)
= duW_{u,v}^{d-1}.\ee Differentiate
both sides of (\ref{fraction3}) with respect to $u$ and $v$ and letting $s\to u^+$ and $t\to v^-$, we arrive at
 \be \label{eqn2} \frac{H'_u}{H_v-H_u}=\frac{1}{v-u}. \ee
Equation (\ref{eqn2}) has the only solution given by
(\ref{exp_type}). If $k=1$ and $r\ge 2$ the proof is similar and
is omitted here. If $k=r=1$, then
(\ref{newstuff2}) simplifies to
\[E\left[
R_n-(d+1)R_m|R_{m-1}=u,R_{n+1}=v\right]=-du. \]
Repeating the arguments for the case $k\ge 2$ above, it is not difficult to obtain equation (\ref{eqn2}) with only solution (\ref{exp_type}).
The sufficiency is proved.

{\bf Necessity.} Using Lemma~2, it is not difficult to verify that the distribution function
 $F(x)=1-\exp\{-c(x-l_F)\}$ satisfies
(\ref{newstuff2}). The theorem
is proved.

\begin{center}
\hspace{0.5cm}
\end{center}

\small{
\begin{tabular}{ll}

\noindent George P. Yanev & M. Ahsanullah \\
Department of Mathematics & Department of Management Sciences \\
University of Texas - Pan American &  Rider University \\
1201 West University Drive & 2083 Lawrenceville Road  \\
Edinburg, TX 78539, USA & Lawrenceville, NJ 08648, USA \\
E-mail: yanevgp@utpa.edu  & E-mail: ahsan@rider.edu
\end{tabular}

}

\end{document}